\documentclass[a4paper, 11pt]{amsart}
\usepackage{epic, eepic,amssymb,epsf}

\input xy
\xyoption{all}

\theoremstyle{definition}

\def\bC{\mathbb{C}}

\def\bN{\mathbb{N}}

\def\cA{\mathcal{A}}

\def\cC{\mathcal{C}}
\def\cS{\mathcal{S}}

\def\cT{\mathcal{T}}

\newcommand{\ra}{\rightarrow}

\newcommand{\hqft}{homotopy quantum field theory }
\newcommand{\hqfts}{homotopy quantum field theories }
\newcommand{\tift}{thin-invariant field theory }
\newcommand{\tifts}{thin-invariant field theories }

\title{Homotopy Quantum Field Theories and Related Ideas}
\author{Mark Brightwell}

\author{Paul Turner}
\address{Heriot-Watt University\\ Edinburgh EH14 4AS\\Scotland}

\author{Simon Willerton}

\begin{document}

\maketitle

\begin{abstract}
In this short note we provide a review of some
developments in the area of homotopy quantum field theories, loosely
based on a talk given by the second author at the Xth Oporto Meeting
on Geometry, Topology and Physics.
\end{abstract}

\section*{Homotopy Quantum Field Theories}
Homotopy quantum field theories were invented by Turaev
\cite{Turaev:HomotopyFieldTheoryInDimension2}, though the idea goes
back to Segal's discussion of the possible geometry underlying
elliptic cohomology \cite{se:ec}. Segal's construction is a
generalisation of his definition of conformal field theory to the
situation where one has a target or background space $X$. He assigns a
topological vector space $E(\gamma)$ to each collection of loops
$\gamma$ in a space $X$ and a trace-class map $E(\sigma)\colon
E(\gamma) \ra E(\gamma^\prime)$ to each Riemann surface $\Sigma$
equipped with a map $\sigma\colon \Sigma \ra X$ agreeing with
$\gamma^{op} \sqcup \gamma^\prime$ on the boundary. The assignment is
multiplicative in the sense that $E(\gamma_1 \sqcup \gamma_2)$ is
isomorphic to $E(\gamma_1 )\otimes E(\gamma_2)$. The result can be
thought of as a kind of infinite dimensional bundle on the free loop
space of $X$, together with a generalised connection which describes
``parallel transport'' along surfaces.

\vspace{0.3cm}
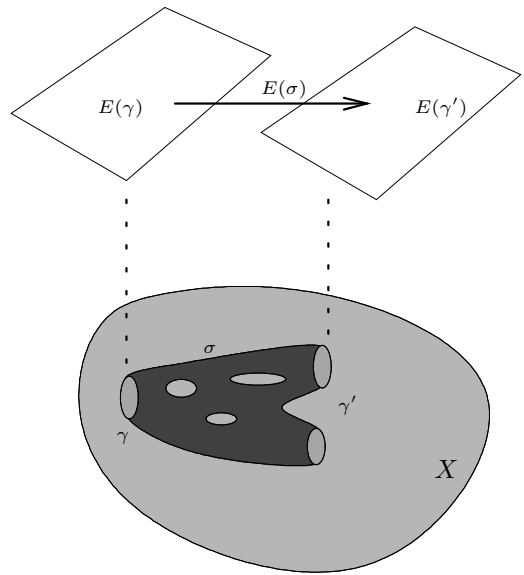
\begin{center}
\setlength{\unitlength}{0.00033333in}
\begingroup\makeatletter\ifx\SetFigFont\undefined%
\gdef\SetFigFont#1#2#3#4#5{%
  \reset@font\fontsize{#1}{#2pt}%
  \fontfamily{#3}\fontseries{#4}\fontshape{#5}%
  \selectfont}%
\fi\endgroup%
{\renewcommand{\dashlinestretch}{30}
\begin{picture}(7974,8855)(0,-10)
\texture{88555555 55000000 555555 55000000 555555 55000000 555555 55000000 
	555555 55000000 555555 55000000 555555 55000000 555555 55000000 
	555555 55000000 555555 55000000 555555 55000000 555555 55000000 
	555555 55000000 555555 55000000 555555 55000000 555555 55000000 }
\shade\path(1988,994)(2026,964)(2065,934)
	(2106,904)(2147,875)(2190,845)
	(2233,815)(2278,786)(2324,756)
	(2372,727)(2420,698)(2470,669)
	(2521,640)(2574,611)(2627,583)
	(2682,555)(2738,527)(2795,499)
	(2853,472)(2912,445)(2972,418)
	(3033,392)(3096,367)(3159,342)
	(3223,318)(3287,294)(3353,271)
	(3419,249)(3485,227)(3552,207)
	(3620,187)(3687,168)(3755,150)
	(3824,133)(3892,117)(3961,102)
	(4029,88)(4097,75)(4165,64)
	(4233,53)(4301,44)(4368,35)
	(4435,28)(4501,23)(4567,18)
	(4632,15)(4697,12)(4761,12)
	(4824,12)(4887,14)(4949,16)
	(5011,21)(5071,26)(5131,33)
	(5191,41)(5249,50)(5307,61)
	(5365,73)(5422,86)(5478,101)
	(5534,117)(5591,135)(5647,155)
	(5704,176)(5760,199)(5816,223)
	(5871,249)(5926,276)(5981,305)
	(6035,335)(6089,367)(6143,400)
	(6196,435)(6249,472)(6302,509)
	(6354,548)(6405,589)(6456,631)
	(6507,674)(6556,719)(6605,765)
	(6653,811)(6701,859)(6747,908)
	(6792,958)(6837,1009)(6880,1061)
	(6922,1114)(6963,1167)(7002,1221)
	(7041,1275)(7078,1330)(7113,1385)
	(7147,1440)(7179,1496)(7210,1551)
	(7239,1607)(7267,1663)(7293,1718)
	(7317,1774)(7340,1829)(7361,1883)
	(7380,1938)(7397,1991)(7413,2045)
	(7426,2098)(7438,2150)(7449,2201)
	(7457,2252)(7464,2302)(7469,2352)
	(7472,2401)(7474,2449)(7474,2496)
	(7472,2543)(7468,2588)(7462,2634)
	(7455,2678)(7446,2722)(7435,2766)
	(7422,2811)(7407,2854)(7391,2897)
	(7372,2939)(7351,2981)(7329,3023)
	(7305,3064)(7278,3105)(7250,3145)
	(7220,3185)(7189,3225)(7155,3264)
	(7120,3303)(7082,3341)(7044,3379)
	(7003,3417)(6961,3454)(6917,3490)
	(6872,3526)(6825,3562)(6776,3597)
	(6727,3632)(6676,3666)(6624,3699)
	(6570,3731)(6516,3763)(6461,3795)
	(6405,3825)(6348,3855)(6291,3884)
	(6232,3913)(6174,3940)(6115,3967)
	(6056,3993)(5996,4018)(5937,4042)
	(5877,4065)(5818,4088)(5758,4110)
	(5699,4131)(5640,4151)(5582,4170)
	(5523,4189)(5466,4206)(5408,4223)
	(5351,4240)(5295,4255)(5239,4270)
	(5184,4284)(5130,4297)(5075,4310)
	(5022,4322)(4969,4334)(4917,4345)
	(4865,4355)(4805,4367)(4745,4378)
	(4686,4388)(4628,4397)(4569,4406)
	(4511,4415)(4453,4422)(4395,4430)
	(4337,4436)(4279,4442)(4222,4448)
	(4164,4453)(4107,4457)(4050,4461)
	(3993,4465)(3936,4468)(3880,4470)
	(3823,4472)(3768,4473)(3712,4474)
	(3658,4475)(3603,4475)(3550,4474)
	(3496,4474)(3444,4472)(3393,4471)
	(3342,4469)(3292,4466)(3243,4464)
	(3196,4461)(3149,4457)(3103,4454)
	(3058,4450)(3015,4446)(2972,4442)
	(2931,4438)(2891,4433)(2852,4429)
	(2813,4424)(2776,4419)(2740,4414)
	(2705,4409)(2671,4404)(2638,4399)
	(2606,4394)(2575,4389)(2544,4384)
	(2514,4379)(2470,4371)(2427,4364)
	(2386,4356)(2345,4349)(2305,4341)
	(2266,4333)(2227,4324)(2189,4316)
	(2151,4307)(2115,4297)(2079,4288)
	(2043,4278)(2008,4267)(1975,4257)
	(1942,4245)(1909,4234)(1878,4222)
	(1848,4210)(1818,4197)(1790,4184)
	(1763,4170)(1736,4156)(1711,4142)
	(1686,4128)(1663,4113)(1640,4098)
	(1618,4082)(1597,4066)(1577,4049)
	(1557,4032)(1537,4015)(1518,3997)
	(1500,3979)(1482,3961)(1464,3941)
	(1446,3921)(1429,3900)(1411,3878)
	(1393,3855)(1376,3831)(1358,3805)
	(1341,3779)(1323,3751)(1306,3722)
	(1289,3692)(1273,3661)(1257,3629)
	(1241,3596)(1226,3562)(1211,3527)
	(1197,3492)(1183,3455)(1170,3418)
	(1158,3381)(1146,3343)(1136,3304)
	(1126,3265)(1117,3226)(1108,3186)
	(1101,3146)(1094,3106)(1088,3065)
	(1083,3024)(1078,2983)(1074,2940)
	(1072,2898)(1070,2866)(1069,2833)
	(1068,2800)(1068,2766)(1068,2731)
	(1069,2696)(1070,2660)(1073,2623)
	(1076,2586)(1079,2547)(1084,2508)
	(1090,2468)(1096,2428)(1104,2386)
	(1112,2344)(1122,2302)(1132,2258)
	(1144,2215)(1157,2170)(1172,2125)
	(1187,2080)(1204,2035)(1222,1989)
	(1241,1943)(1261,1897)(1283,1851)
	(1306,1805)(1331,1759)(1356,1713)
	(1383,1668)(1412,1622)(1441,1577)
	(1472,1533)(1504,1489)(1537,1445)
	(1571,1402)(1607,1359)(1644,1316)
	(1682,1275)(1722,1233)(1762,1192)
	(1805,1152)(1848,1112)(1893,1072)
	(1940,1033)(1988,994)
\path(1988,994)(2026,964)(2065,934)
	(2106,904)(2147,875)(2190,845)
	(2233,815)(2278,786)(2324,756)
	(2372,727)(2420,698)(2470,669)
	(2521,640)(2574,611)(2627,583)
	(2682,555)(2738,527)(2795,499)
	(2853,472)(2912,445)(2972,418)
	(3033,392)(3096,367)(3159,342)
	(3223,318)(3287,294)(3353,271)
	(3419,249)(3485,227)(3552,207)
	(3620,187)(3687,168)(3755,150)
	(3824,133)(3892,117)(3961,102)
	(4029,88)(4097,75)(4165,64)
	(4233,53)(4301,44)(4368,35)
	(4435,28)(4501,23)(4567,18)
	(4632,15)(4697,12)(4761,12)
	(4824,12)(4887,14)(4949,16)
	(5011,21)(5071,26)(5131,33)
	(5191,41)(5249,50)(5307,61)
	(5365,73)(5422,86)(5478,101)
	(5534,117)(5591,135)(5647,155)
	(5704,176)(5760,199)(5816,223)
	(5871,249)(5926,276)(5981,305)
	(6035,335)(6089,367)(6143,400)
	(6196,435)(6249,472)(6302,509)
	(6354,548)(6405,589)(6456,631)
	(6507,674)(6556,719)(6605,765)
	(6653,811)(6701,859)(6747,908)
	(6792,958)(6837,1009)(6880,1061)
	(6922,1114)(6963,1167)(7002,1221)
	(7041,1275)(7078,1330)(7113,1385)
	(7147,1440)(7179,1496)(7210,1551)
	(7239,1607)(7267,1663)(7293,1718)
	(7317,1774)(7340,1829)(7361,1883)
	(7380,1938)(7397,1991)(7413,2045)
	(7426,2098)(7438,2150)(7449,2201)
	(7457,2252)(7464,2302)(7469,2352)
	(7472,2401)(7474,2449)(7474,2496)
	(7472,2543)(7468,2588)(7462,2634)
	(7455,2678)(7446,2722)(7435,2766)
	(7422,2811)(7407,2854)(7391,2897)
	(7372,2939)(7351,2981)(7329,3023)
	(7305,3064)(7278,3105)(7250,3145)
	(7220,3185)(7189,3225)(7155,3264)
	(7120,3303)(7082,3341)(7044,3379)
	(7003,3417)(6961,3454)(6917,3490)
	(6872,3526)(6825,3562)(6776,3597)
	(6727,3632)(6676,3666)(6624,3699)
	(6570,3731)(6516,3763)(6461,3795)
	(6405,3825)(6348,3855)(6291,3884)
	(6232,3913)(6174,3940)(6115,3967)
	(6056,3993)(5996,4018)(5937,4042)
	(5877,4065)(5818,4088)(5758,4110)
	(5699,4131)(5640,4151)(5582,4170)
	(5523,4189)(5466,4206)(5408,4223)
	(5351,4240)(5295,4255)(5239,4270)
	(5184,4284)(5130,4297)(5075,4310)
	(5022,4322)(4969,4334)(4917,4345)
	(4865,4355)(4805,4367)(4745,4378)
	(4686,4388)(4628,4397)(4569,4406)
	(4511,4415)(4453,4422)(4395,4430)
	(4337,4436)(4279,4442)(4222,4448)
	(4164,4453)(4107,4457)(4050,4461)
	(3993,4465)(3936,4468)(3880,4470)
	(3823,4472)(3768,4473)(3712,4474)
	(3658,4475)(3603,4475)(3550,4474)
	(3496,4474)(3444,4472)(3393,4471)
	(3342,4469)(3292,4466)(3243,4464)
	(3196,4461)(3149,4457)(3103,4454)
	(3058,4450)(3015,4446)(2972,4442)
	(2931,4438)(2891,4433)(2852,4429)
	(2813,4424)(2776,4419)(2740,4414)
	(2705,4409)(2671,4404)(2638,4399)
	(2606,4394)(2575,4389)(2544,4384)
	(2514,4379)(2470,4371)(2427,4364)
	(2386,4356)(2345,4349)(2305,4341)
	(2266,4333)(2227,4324)(2189,4316)
	(2151,4307)(2115,4297)(2079,4288)
	(2043,4278)(2008,4267)(1975,4257)
	(1942,4245)(1909,4234)(1878,4222)
	(1848,4210)(1818,4197)(1790,4184)
	(1763,4170)(1736,4156)(1711,4142)
	(1686,4128)(1663,4113)(1640,4098)
	(1618,4082)(1597,4066)(1577,4049)
	(1557,4032)(1537,4015)(1518,3997)
	(1500,3979)(1482,3961)(1464,3941)
	(1446,3921)(1429,3900)(1411,3878)
	(1393,3855)(1376,3831)(1358,3805)
	(1341,3779)(1323,3751)(1306,3722)
	(1289,3692)(1273,3661)(1257,3629)
	(1241,3596)(1226,3562)(1211,3527)
	(1197,3492)(1183,3455)(1170,3418)
	(1158,3381)(1146,3343)(1136,3304)
	(1126,3265)(1117,3226)(1108,3186)
	(1101,3146)(1094,3106)(1088,3065)
	(1083,3024)(1078,2983)(1074,2940)
	(1072,2898)(1070,2866)(1069,2833)
	(1068,2800)(1068,2766)(1068,2731)
	(1069,2696)(1070,2660)(1073,2623)
	(1076,2586)(1079,2547)(1084,2508)
	(1090,2468)(1096,2428)(1104,2386)
	(1112,2344)(1122,2302)(1132,2258)
	(1144,2215)(1157,2170)(1172,2125)
	(1187,2080)(1204,2035)(1222,1989)
	(1241,1943)(1261,1897)(1283,1851)
	(1306,1805)(1331,1759)(1356,1713)
	(1383,1668)(1412,1622)(1441,1577)
	(1472,1533)(1504,1489)(1537,1445)
	(1571,1402)(1607,1359)(1644,1316)
	(1682,1275)(1722,1233)(1762,1192)
	(1805,1152)(1848,1112)(1893,1072)
	(1940,1033)(1988,994)
\texture{aaffffff ffaaaaaa aaffffff ffaaaaaa aaffffff ffaaaaaa aaffffff ffaaaaaa 
	aaffffff ffaaaaaa aaffffff ffaaaaaa aaffffff ffaaaaaa aaffffff ffaaaaaa 
	aaffffff ffaaaaaa aaffffff ffaaaaaa aaffffff ffaaaaaa aaffffff ffaaaaaa 
	aaffffff ffaaaaaa aaffffff ffaaaaaa aaffffff ffaaaaaa aaffffff ffaaaaaa }
\shade\path(4651,3559)(4684,3560)(4715,3559)
	(4743,3558)(4769,3556)(4793,3553)
	(4815,3549)(4835,3545)(4853,3540)
	(4870,3534)(4885,3527)(4898,3520)
	(4909,3512)(4919,3503)(4927,3493)
	(4934,3484)(4940,3473)(4944,3462)
	(4947,3451)(4950,3439)(4951,3428)
	(4952,3415)(4952,3403)(4952,3391)
	(4951,3379)(4950,3366)(4949,3354)
	(4948,3341)(4947,3329)(4947,3316)
	(4946,3304)(4947,3287)(4948,3270)
	(4950,3252)(4952,3234)(4954,3216)
	(4956,3197)(4959,3178)(4960,3158)
	(4962,3139)(4962,3119)(4962,3099)
	(4960,3079)(4957,3059)(4953,3040)
	(4947,3021)(4939,3003)(4930,2985)
	(4918,2968)(4905,2952)(4889,2936)
	(4871,2920)(4850,2905)(4831,2892)
	(4809,2880)(4786,2867)(4761,2855)
	(4734,2843)(4705,2830)(4675,2817)
	(4645,2804)(4613,2792)(4580,2779)
	(4548,2766)(4515,2753)(4483,2740)
	(4452,2727)(4422,2714)(4393,2702)
	(4366,2689)(4341,2677)(4318,2665)
	(4298,2653)(4281,2642)(4266,2630)
	(4254,2619)(4245,2608)(4239,2597)
	(4237,2586)(4238,2574)(4243,2562)
	(4251,2550)(4263,2537)(4278,2524)
	(4296,2511)(4318,2498)(4342,2484)
	(4368,2470)(4396,2455)(4427,2441)
	(4458,2426)(4490,2411)(4523,2397)
	(4556,2382)(4588,2367)(4620,2353)
	(4650,2339)(4679,2325)(4706,2312)
	(4732,2299)(4755,2286)(4776,2273)
	(4794,2260)(4813,2245)(4829,2229)
	(4842,2213)(4853,2197)(4861,2181)
	(4866,2164)(4870,2146)(4871,2129)
	(4871,2111)(4869,2093)(4866,2075)
	(4862,2057)(4858,2040)(4853,2022)
	(4848,2005)(4844,1989)(4840,1973)
	(4837,1957)(4835,1941)(4834,1926)
	(4835,1916)(4836,1905)(4837,1895)
	(4838,1885)(4840,1874)(4842,1864)
	(4844,1853)(4846,1843)(4848,1832)
	(4849,1821)(4850,1811)(4850,1801)
	(4849,1790)(4847,1780)(4845,1771)
	(4841,1761)(4836,1752)(4829,1743)
	(4821,1735)(4812,1727)(4801,1720)
	(4788,1713)(4774,1707)(4758,1701)
	(4740,1696)(4720,1692)(4698,1688)
	(4674,1684)(4648,1681)(4619,1679)
	(4596,1677)(4571,1676)(4545,1675)
	(4517,1674)(4488,1674)(4457,1674)
	(4424,1674)(4390,1675)(4354,1676)
	(4317,1677)(4278,1679)(4238,1681)
	(4197,1683)(4154,1686)(4110,1689)
	(4065,1693)(4018,1696)(3971,1700)
	(3924,1705)(3875,1710)(3826,1715)
	(3777,1720)(3728,1726)(3678,1732)
	(3629,1738)(3579,1745)(3530,1752)
	(3482,1759)(3434,1767)(3386,1774)
	(3339,1782)(3293,1790)(3248,1798)
	(3203,1807)(3160,1815)(3117,1824)
	(3075,1833)(3034,1842)(2993,1852)
	(2953,1862)(2907,1873)(2862,1886)
	(2817,1898)(2773,1911)(2730,1925)
	(2686,1939)(2643,1953)(2601,1968)
	(2559,1984)(2518,2000)(2477,2016)
	(2436,2033)(2397,2050)(2358,2067)
	(2320,2085)(2283,2103)(2247,2121)
	(2212,2139)(2179,2157)(2147,2175)
	(2116,2193)(2086,2211)(2058,2229)
	(2032,2246)(2007,2263)(1984,2280)
	(1962,2297)(1941,2313)(1922,2329)
	(1904,2344)(1888,2360)(1873,2374)
	(1859,2389)(1846,2403)(1828,2425)
	(1813,2446)(1800,2468)(1788,2489)
	(1779,2510)(1771,2531)(1766,2552)
	(1761,2572)(1758,2593)(1757,2614)
	(1757,2634)(1757,2653)(1759,2673)
	(1761,2691)(1764,2709)(1767,2727)
	(1770,2743)(1773,2759)(1776,2775)
	(1779,2789)(1781,2804)(1782,2818)
	(1783,2835)(1783,2852)(1782,2869)
	(1780,2886)(1779,2903)(1777,2920)
	(1776,2937)(1775,2955)(1775,2971)
	(1776,2988)(1779,3004)(1784,3019)
	(1790,3034)(1798,3048)(1808,3061)
	(1821,3073)(1836,3085)(1854,3096)
	(1870,3104)(1887,3112)(1906,3120)
	(1927,3128)(1950,3135)(1974,3143)
	(2000,3150)(2027,3157)(2055,3165)
	(2085,3172)(2115,3178)(2145,3185)
	(2176,3191)(2207,3197)(2237,3203)
	(2267,3208)(2295,3213)(2323,3218)
	(2349,3222)(2375,3226)(2398,3230)
	(2420,3233)(2441,3236)(2460,3239)
	(2477,3242)(2494,3244)(2509,3247)
	(2524,3249)(2539,3251)(2554,3253)
	(2569,3256)(2584,3258)(2600,3261)
	(2617,3264)(2635,3266)(2654,3270)
	(2675,3273)(2697,3277)(2721,3281)
	(2747,3285)(2774,3290)(2804,3295)
	(2835,3301)(2870,3307)(2906,3313)
	(2945,3320)(2988,3327)(3033,3335)
	(3063,3340)(3093,3345)(3125,3351)
	(3159,3357)(3194,3363)(3230,3369)
	(3268,3375)(3307,3382)(3347,3389)
	(3389,3396)(3432,3403)(3476,3410)
	(3522,3418)(3568,3425)(3615,3433)
	(3662,3440)(3710,3448)(3759,3455)
	(3808,3463)(3857,3470)(3905,3477)
	(3954,3484)(4002,3491)(4050,3498)
	(4096,3504)(4143,3511)(4188,3517)
	(4232,3522)(4275,3527)(4316,3532)
	(4356,3537)(4395,3541)(4433,3545)
	(4469,3548)(4503,3551)(4536,3553)
	(4567,3555)(4596,3557)(4624,3558)(4651,3559)
\path(4651,3559)(4684,3560)(4715,3559)
	(4743,3558)(4769,3556)(4793,3553)
	(4815,3549)(4835,3545)(4853,3540)
	(4870,3534)(4885,3527)(4898,3520)
	(4909,3512)(4919,3503)(4927,3493)
	(4934,3484)(4940,3473)(4944,3462)
	(4947,3451)(4950,3439)(4951,3428)
	(4952,3415)(4952,3403)(4952,3391)
	(4951,3379)(4950,3366)(4949,3354)
	(4948,3341)(4947,3329)(4947,3316)
	(4946,3304)(4947,3287)(4948,3270)
	(4950,3252)(4952,3234)(4954,3216)
	(4956,3197)(4959,3178)(4960,3158)
	(4962,3139)(4962,3119)(4962,3099)
	(4960,3079)(4957,3059)(4953,3040)
	(4947,3021)(4939,3003)(4930,2985)
	(4918,2968)(4905,2952)(4889,2936)
	(4871,2920)(4850,2905)(4831,2892)
	(4809,2880)(4786,2867)(4761,2855)
	(4734,2843)(4705,2830)(4675,2817)
	(4645,2804)(4613,2792)(4580,2779)
	(4548,2766)(4515,2753)(4483,2740)
	(4452,2727)(4422,2714)(4393,2702)
	(4366,2689)(4341,2677)(4318,2665)
	(4298,2653)(4281,2642)(4266,2630)
	(4254,2619)(4245,2608)(4239,2597)
	(4237,2586)(4238,2574)(4243,2562)
	(4251,2550)(4263,2537)(4278,2524)
	(4296,2511)(4318,2498)(4342,2484)
	(4368,2470)(4396,2455)(4427,2441)
	(4458,2426)(4490,2411)(4523,2397)
	(4556,2382)(4588,2367)(4620,2353)
	(4650,2339)(4679,2325)(4706,2312)
	(4732,2299)(4755,2286)(4776,2273)
	(4794,2260)(4813,2245)(4829,2229)
	(4842,2213)(4853,2197)(4861,2181)
	(4866,2164)(4870,2146)(4871,2129)
	(4871,2111)(4869,2093)(4866,2075)
	(4862,2057)(4858,2040)(4853,2022)
	(4848,2005)(4844,1989)(4840,1973)
	(4837,1957)(4835,1941)(4834,1926)
	(4835,1916)(4836,1905)(4837,1895)
	(4838,1885)(4840,1874)(4842,1864)
	(4844,1853)(4846,1843)(4848,1832)
	(4849,1821)(4850,1811)(4850,1801)
	(4849,1790)(4847,1780)(4845,1771)
	(4841,1761)(4836,1752)(4829,1743)
	(4821,1735)(4812,1727)(4801,1720)
	(4788,1713)(4774,1707)(4758,1701)
	(4740,1696)(4720,1692)(4698,1688)
	(4674,1684)(4648,1681)(4619,1679)
	(4596,1677)(4571,1676)(4545,1675)
	(4517,1674)(4488,1674)(4457,1674)
	(4424,1674)(4390,1675)(4354,1676)
	(4317,1677)(4278,1679)(4238,1681)
	(4197,1683)(4154,1686)(4110,1689)
	(4065,1693)(4018,1696)(3971,1700)
	(3924,1705)(3875,1710)(3826,1715)
	(3777,1720)(3728,1726)(3678,1732)
	(3629,1738)(3579,1745)(3530,1752)
	(3482,1759)(3434,1767)(3386,1774)
	(3339,1782)(3293,1790)(3248,1798)
	(3203,1807)(3160,1815)(3117,1824)
	(3075,1833)(3034,1842)(2993,1852)
	(2953,1862)(2907,1873)(2862,1886)
	(2817,1898)(2773,1911)(2730,1925)
	(2686,1939)(2643,1953)(2601,1968)
	(2559,1984)(2518,2000)(2477,2016)
	(2436,2033)(2397,2050)(2358,2067)
	(2320,2085)(2283,2103)(2247,2121)
	(2212,2139)(2179,2157)(2147,2175)
	(2116,2193)(2086,2211)(2058,2229)
	(2032,2246)(2007,2263)(1984,2280)
	(1962,2297)(1941,2313)(1922,2329)
	(1904,2344)(1888,2360)(1873,2374)
	(1859,2389)(1846,2403)(1828,2425)
	(1813,2446)(1800,2468)(1788,2489)
	(1779,2510)(1771,2531)(1766,2552)
	(1761,2572)(1758,2593)(1757,2614)
	(1757,2634)(1757,2653)(1759,2673)
	(1761,2691)(1764,2709)(1767,2727)
	(1770,2743)(1773,2759)(1776,2775)
	(1779,2789)(1781,2804)(1782,2818)
	(1783,2835)(1783,2852)(1782,2869)
	(1780,2886)(1779,2903)(1777,2920)
	(1776,2937)(1775,2955)(1775,2971)
	(1776,2988)(1779,3004)(1784,3019)
	(1790,3034)(1798,3048)(1808,3061)
	(1821,3073)(1836,3085)(1854,3096)
	(1870,3104)(1887,3112)(1906,3120)
	(1927,3128)(1950,3135)(1974,3143)
	(2000,3150)(2027,3157)(2055,3165)
	(2085,3172)(2115,3178)(2145,3185)
	(2176,3191)(2207,3197)(2237,3203)
	(2267,3208)(2295,3213)(2323,3218)
	(2349,3222)(2375,3226)(2398,3230)
	(2420,3233)(2441,3236)(2460,3239)
	(2477,3242)(2494,3244)(2509,3247)
	(2524,3249)(2539,3251)(2554,3253)
	(2569,3256)(2584,3258)(2600,3261)
	(2617,3264)(2635,3266)(2654,3270)
	(2675,3273)(2697,3277)(2721,3281)
	(2747,3285)(2774,3290)(2804,3295)
	(2835,3301)(2870,3307)(2906,3313)
	(2945,3320)(2988,3327)(3033,3335)
	(3063,3340)(3093,3345)(3125,3351)
	(3159,3357)(3194,3363)(3230,3369)
	(3268,3375)(3307,3382)(3347,3389)
	(3389,3396)(3432,3403)(3476,3410)
	(3522,3418)(3568,3425)(3615,3433)
	(3662,3440)(3710,3448)(3759,3455)
	(3808,3463)(3857,3470)(3905,3477)
	(3954,3484)(4002,3491)(4050,3498)
	(4096,3504)(4143,3511)(4188,3517)
	(4232,3522)(4275,3527)(4316,3532)
	(4356,3537)(4395,3541)(4433,3545)
	(4469,3548)(4503,3551)(4536,3553)
	(4567,3555)(4596,3557)(4624,3558)(4651,3559)
\texture{80222222 22555555 55808080 80555555 55222222 22555555 55880888 8555555 
	55222222 22555555 55808080 80555555 55222222 22555555 55080808 8555555 
	55222222 22555555 55808080 80555555 55222222 22555555 55880888 8555555 
	55222222 22555555 55808080 80555555 55222222 22555555 55080808 8555555 }
\put(4771,1973){\shade\ellipse{286}{574}}
\put(4771,1973){\ellipse{286}{574}}
\texture{88555555 55000000 555555 55000000 555555 55000000 555555 55000000 
	555555 55000000 555555 55000000 555555 55000000 555555 55000000 
	555555 55000000 555555 55000000 555555 55000000 555555 55000000 
	555555 55000000 555555 55000000 555555 55000000 555555 55000000 }
\put(2667,2881){\shade\ellipse{478}{286}}
\put(2667,2881){\ellipse{478}{286}}
\put(3289,2403){\shade\ellipse{478}{190}}
\put(3289,2403){\ellipse{478}{190}}
\put(3863,3025){\shade\ellipse{862}{190}}
\put(3863,3025){\ellipse{862}{190}}
\thicklines
\dashline{90.000}(1812,3278)(1812,5828)
\dashline{90.000}(4962,3728)(4962,5828)
\thinlines
\path(12,7178)(2412,8828)(4062,8078)
	(1812,6128)(12,7178)
\path(3912,6878)(6312,8528)(7962,7778)
	(5712,5828)(3912,6878)
\thicklines
\path(2562,7328)(5562,7328)
\path(5202.000,7238.000)(5562.000,7328.000)(5202.000,7418.000)
\put(1362,7178){\makebox(0,0)[lb]{\smash{{{\SetFigFont{7}{8.4}{\rmdefault}{\mddefault}{\updefault}$E(\gamma)$}}}}}
\put(6312,7178){\makebox(0,0)[lb]{\smash{{{\SetFigFont{7}{8.4}{\rmdefault}{\mddefault}{\updefault}$E(\gamma^\prime)$}}}}}
\put(3912,7478){\makebox(0,0)[lb]{\smash{{{\SetFigFont{7}{8.4}{\rmdefault}{\mddefault}{\updefault}$E(\sigma)$}}}}}
\put(1662,2078){\makebox(0,0)[lb]{\smash{{{\SetFigFont{7}{8.4}{\rmdefault}{\mddefault}{\updefault}$\gamma$}}}}}
\put(5112,2528){\makebox(0,0)[lb]{\smash{{{\SetFigFont{7}{8.4}{\rmdefault}{\mddefault}{\updefault}$\gamma^\prime$}}}}}
\put(3012,3428){\makebox(0,0)[lb]{\smash{{{\SetFigFont{7}{8.4}{\rmdefault}{\mddefault}{\updefault}$\sigma$}}}}}
\put(6612,1478){\makebox(0,0)[lb]{\smash{{{\SetFigFont{10}{12.0}{\rmdefault}{\mddefault}{\updefault}$X$}}}}}
\thinlines
\texture{80222222 22555555 55808080 80555555 55222222 22555555 55880888 8555555 
	55222222 22555555 55808080 80555555 55222222 22555555 55080808 8555555 
	55222222 22555555 55808080 80555555 55222222 22555555 55880888 8555555 
	55222222 22555555 55808080 80555555 55222222 22555555 55080808 8555555 }
\put(4867,3216){\shade\ellipse{286}{668}}
\put(4867,3216){\ellipse{286}{668}}
\put(1854,2738){\shade\ellipse{286}{668}}
\put(1854,2738){\ellipse{286}{668}}
\end{picture}
}

\end{center}
\vspace{0.3cm}

One can package this in terms of representations of a category $\cC_X$
whose objects are pairs $(\Gamma, \gamma)$, where $\Gamma$ is a
compact, closed, oriented 1-manifold and $\gamma\colon \Gamma \ra X$
is a continuous function; and whose morphisms are equivalence classes
of triples $(\Sigma, \alpha, \sigma)$, where $\Sigma$ is a Riemann
surface, $\alpha$ a boundary identification
$\alpha\colon\partial\Sigma \cong \Gamma^{op} \sqcup \Gamma^\prime$
and $\sigma\colon \Sigma \ra X$ is a continuous function equal to
$\gamma$ and $ \gamma^\prime$ on the boundary.  Segal defines an
elliptic object to be a multiplicative functor from this category to
the category of topological vector spaces satisfying a number of
further conditions.  For $X$ a point one regains the category at the
heart of Segal's axiomatic definition of conformal field theory.  From one
point of view the above is a way of organising the surfaces and maps
to an auxiliary space $X$ which occur in non-linear $\sigma$-models in
string theory, but stopping short of integrating over mapping
spaces.

A 1+1-dimensional {\em \hqft } is a variant of the above, in which the
complex structure is neglected, topological vector spaces are replaced
with finite dimensional complex vector spaces and linear maps
associated to cobordisms are invariant under deformation by homotopy
(relative to the boundary). These can be thought of geometrically as
flat ``higher'' bundles with base the free loop space. In general an
$n+1$-dimensional \hqft is an $n+1$-dimensional topological quantum
field theory ``with background''. Of particular interest is when the
background space is a classifying space of some kind and so maps to it
have a geometrical meaning, for example $X=BG$ in which case maps to
$X$ are interpreted as $G$-bundles. For $G$ a finite group this makes
contact with Dijkgraaf-Witten theory.

For a path connected pointed space $X$, define an $X$-manifold to be a
closed oriented $n$-manifold $M$ with pointed components, equipped
with a based map $\gamma\colon M\ra X$. Define an $X$-homeomorphism of
$X$-manifolds to be an orientation preserving homeomorphism $f\colon
M\ra M^\prime$ sending basepoints to basepoints such that
$\gamma^\prime \circ f = \gamma$.  An $X$-cobordism is an oriented
$n+1$-manifold whose boundary is an $X$-manifold (where the
orientation of the ingoing boundary components is opposite to the
induced one) together with a map to $X$ (taking boundary basepoints to
the basepoint in $X$). Where boundaries agree, $X$-cobordisms can be
glued using $X$-homeomorphisms. Let $W_0 \cup_f W_1$ denote the result
of gluing $W_0$ to $W_1$ using the $X$-homeomorphism $f$ where the
outgoing boundary of $W_0$ is glued to the ingoing boundary of $W_1$.
Turaev then defines an $n+1$-dimensional \hqft (over $\bC$) with
target $X$ as an assignment of a finite dimensional vector space $A_M$
to each $X$-manifold $M$, of an isomorphism $f_\sharp\colon M\ra
M^\prime$ to each $X$-homeomorphism $f\colon M\ra M^\prime$, and of a
linear map $\tau(W)\colon A_{M_0} \ra A_{M_1}$ to each $X$-cobordism
$W$ from $M_0$ to $M_1$, satisfying the following axioms, where the
numbering is Turaev's.
\begin{itemize}
\item[(1.2.1)] each isomorphism $f_\sharp$ is invariant under
  isotopies of $X$-homeomorphisms and $(f^\prime f)_\sharp =
  f^\prime_\sharp f_\sharp$
\item[(1.2.2)] for two $X$-manifolds $M$ and $N$ there is a symmetric
  natural isomorphism $A_{M\sqcup N} \cong A_M \otimes A_N$
\item[(1.2.3)] there is an isomorphism $A_\emptyset \cong \bC$
\item[(1.2.4)] the maps $\tau$ are natural with respect to
  $X$-homeomorphisms
\item[(1.2.5)] for two $X$-cobordisms $W_1$ and $W_2$, $\tau
  (W_1\sqcup W_2) = \tau (W_1) \otimes \tau(W_2)$
\item[(1.2.6)] $\tau (W_0 \cup_f W_1) = \tau(W_1)\circ f_\sharp \circ
  \tau(W_0)$ 
\item[(1.2.7)] for any $X$-manifold $\gamma\colon M\ra X$, the
  $X$-cobordism $\sigma\colon M\times [0,1]\ra X$ given by $(x,t)\mapsto
  \gamma(x)$ is assigned the identity map $id\colon A_M \ra A_M$.
\item[(1.2.8)] for any $X$-cobordism $W$, $\tau(W)$ is invariant under
  homotopies relative to $\partial W$.
\end{itemize}

In the above there is one important difference from Turaev's
definition as stated in \cite{Turaev:HomotopyFieldTheoryInDimension2},
namely part of his axiom (1.2.7) has been omitted. In his manuscript
he demands $\tau (M\times [0,1]\ra X) = id$ for {\em any}
$X$-cobordism of the form $M\times [0,1]\ra X$ which by force makes
homotopy information in dimensions greater than $n$ redundant. By
omitting this (and thus giving a role to higher homotopy)
Turaev's results in \cite{Turaev:HomotopyFieldTheoryInDimension2} must
be restated as applying to $X$ an Eilenberg-Maclane space only. This
is the position adopted by Rodrigues in \cite{ro:hqft}, where a
careful discussion of the axioms can be found.

Notice that for a contractible target space the above reduces to a
topological quantum field theory as for each manifold or cobordism
there is a unique up to homotopy map to $X$. Notice too that the
definition provides numerical invariants of closed $X$-cobordisms by
methods standard in TQFT: regard a closed manifold as a cobordism from
the empty manifold to itself which gives a linear map $\bC \cong
A_\emptyset \ra A_\emptyset \cong \bC$ whose value on 1 is the desired
invariant.  By methods standard in TQFT it can also be seen that
$A_{M^{op}} \cong A_M^*$. Also, associated to each point $x$ in $X$ there is an
induced TQFT defined by those maps collapsing all manifolds and
cobordisms to $x$.  

Concurrent with Turaev's definition, the first two authors defined
a categorical version of the above in dimension $1+1$.  Letting $S_n$
denote $n$ copies of a standard circle define the {\em homotopy
  surface category
$\cS_X$ of a space $X$} to be the category with the following objects
and morphisms.
\begin{itemize}
\item Objects are pairs $(n,s)$ where $n\in \bN$  
and $s\colon S_n\ra X$ is a 
continuous function.
\item Morphisms from $(n,s)$ to $(n^\prime,s^\prime)$ are triples $(\Sigma,
  \alpha, \sigma)$ where
\begin{enumerate}
\item $\Sigma$ is an smooth oriented surface
\item $\alpha\colon \partial \Sigma \ra S_n^{op} \sqcup S_{n^\prime}$ is an
  orientation preserving homeomorphism and
\item $\sigma\colon \Sigma \ra X$ is a continuous function such that 
$\sigma|_{\partial \Sigma}\circ \alpha^{-1} = s \sqcup s^\prime$.
\end{enumerate}
A morphism $(\Sigma_1, \alpha_1, \sigma_1)$ is identified with $(\Sigma_2,
\alpha_2, \sigma_2 )$ if there exists a diffeomorphism
$T\colon \Sigma_1 \ra \Sigma_2$ such that  $\alpha_2 \circ
T|_{\partial\Sigma} = \alpha_1$ and the diagram
\begin{equation}\label{diag:morph}
\xymatrix{\Sigma_1 \ar[rr]^{T} \ar[dr]_{g_1} & & \Sigma_2
  \ar[dl]^{g_2}\\ & X}
\end{equation}
commutes up to homotopy relative to
the boundary of $\Sigma_1$.
\end{itemize}

This is a monoidal category under disjoint union and an equivalent
definition of \hqft is as a symmetric monoidal functor from $\cS_X$
to the category of finite dimensional complex vector spaces (monoidal
under tensor product).  Notice that $\cS_X$ plays the role of a
``higher'' fundamental groupoid, and just as representations of the
fundamental groupoid correspond to flat bundles, representations of
$\cS_X$ give flat ``higher'' bundles.

Though more care is needed for higher dimensional theories, Rodrigues
\cite{ro:hqft} has succeeded in producing a
cobordism category for each dimension whose multiplicative
representations give homotopy quantum field theories.

\section*{Examples}
\subsection*{Flat vector bundles} An example of a 0+1-dimensional \hqft is
given by a vector bundle with flat connection: the vector space $E(x)$
associated to a point $x$ in $X$ is the fibre over $x$ and a map $E(x)
\ra E(y)$ associated to a path from $x$ to $y$ is given by parallel
transport.

\subsection*{Examples from cocycles} Turaev
\cite{Turaev:HomotopyFieldTheoryInDimension2} constructs
$1+1$-dimensional \hqfts from $2$-cocycles in $X$. The construction
proceeds as follows. Let $\theta\in
C^2X$ be a 2-cocycle and for a loop $\gamma$ set
\[
A_\gamma = \mbox{Span}_\bC \{a\in C_1(S^1)\; | \;
\mbox{$a$ represents the fundamental class of $S^1$} \} / \sim
\]
where $a \sim b $ if $a=\gamma^*\theta(e)b$ for $e\in C_2(S^1) $ such
that $ \partial e=a-b$.  For a cobordism $g\colon \Sigma\to X$ from
$\gamma_1$ to $\gamma_2$ pick a singular two-cycle representative
$f\in C_2(\Sigma)$ of the fundamental class $[\Sigma]\in
H_2(\Sigma,\partial \Sigma)$ and define $A_{\gamma_1} \ra
A_{\gamma_2}$ by $a\mapsto g^*\theta(f)a^\prime$ where $\partial f =
-a+a^\prime$. If two cocycles differ by a coboundary then the theories
constructed above are isomorphic. This construction produces a {\em
  rank one} theory i.e. each vector space $A_\gamma$ is one
dimensional.

In fact this construction is more general and in a similar way Turaev
constructs $n+1$-dimensional \hqfts starting from $n+1$-cocyles in $X$.

\subsection*{Examples from variants of Frobenius algebras} Further
examples in dimension 1+1 can be obtained from Frobenius algebras with
additional structure. For $G$ a discrete group the first two authors
define a {\em $G$-Frobenius algebra} to be a finite dimensional
commutative Frobenius algebra with $G$-action satisfying $g(ab) =
(ga)b= a(gb)$.  Then for $X$ a simply connected space, a
$\pi_2(X)$-Frobenius algebra $V$ gives a 1+1-dimensional \hqft by
setting $A_\gamma=V$ for a single loop $\gamma$. By choosing
contractions for each loop, all linear maps assigned to $X$-cobordisms
are determined from those cobordisms whose boundary loops are trivial
and these in turn are determined by the $\pi_2(X)$-Frobenius algebra.
In particular there is a canonical one-to-one correspondence between
homotopy classes of maps from a cylinder to $X$ with ends mapped to
the basepoint of $X$ and the group $\pi_2(X)$. The linear map assigned
to the cylinder corresponding to $g\in\pi_2(x)$ is given by
multiplication by $g$. In this case the surface category is equivalent
to a labelled version of the category of surfaces without background,
as found in the definition of 1+1-dimensional TQFT. The morphisms are
labelled by elements of $\pi_2(X)$ and when composing morphisms the
labels add.

Similarly, for a discrete group $\pi$, Turaev has defined a {\em
  crossed $\pi$-algebra} to be a $\pi$-graded algebra
$V=\bigoplus_{\alpha\in \pi} V_\alpha$ together with a bilinear form
$\eta\colon V\otimes V \ra \bC$ and a $\pi$-action $\varphi\colon \pi
\ra \mbox{Aut}(V)$ satisfying
\begin{itemize}
\item $\eta |_{V_\alpha\otimes V_\beta}$ is non-degenerate for $\beta =
  \alpha^{-1}$ and zero otherwise and $\eta (ab,c) = \eta(a,bc)$
\item $\varphi (\beta)$ is an algebra automorphism preserving $\eta$
  such that $\varphi(\beta)(V_\alpha) \subset
  V_{\beta\alpha\beta^{-1}}$ and $\varphi(\beta)|_{V_\beta} = id$
 and $(\varphi(\beta)(a)b = ba$
\item $\mbox{Trace}(c\varphi(\beta)\colon V_\alpha \ra V_\alpha) =
  \mbox{Trace}(\varphi(\alpha^{-1})c\colon V_\beta \ra V_\beta)$.
\end{itemize}
Then for $X=K(\pi,1)$ a crossed $\pi$-algebra $V$ gives a
1+1-dimensional \hqft by setting $A_\gamma = V_{[\gamma]}$ where we
identify homotopy classes of maps from $S^1$ to $K(\pi,1)$ with $\pi$
and $[\gamma]$ denotes the homotopy class of $\gamma$.

In fact more is true and Turaev
\cite{Turaev:HomotopyFieldTheoryInDimension2} has shown that the
category of 1+1-dimensional \hqfts with target $K(\pi,1)$ is
equivalent to the category of crossed $\pi$-algebras. Similarly the
first two authors
\cite{BrightwellTurner:RepresentationOfHomotopySurfaceCategory} have
shown that for a simply connected space $X$ the category of
1+1-dimensional \hqfts with target $X$ is equivalent to the category
of $\pi_2(X)$-Frobenius algebras. Rodrigues \cite{ro:hqft} has reformulated
the latter to state that the homotopy surface category of a simply
connected space is universal for $G$-Frobenius objects. This
formalises the view that surface and diagram categories
encode the axioms for algebraic structures.

\subsection*{State sum examples} Further examples of 1+1-dimensional
theories can be obtained by extending state-sum TQFTs to the \hqft
setting. These lattice models have been defined by Turaev
\cite{Turaev:HomotopyFieldTheoryInDimension2} and studied by him for
$X$ an Eilenberg-Maclane space and studied by Rodrigues \cite{ro:hqft}
for simply connected background. In this approach one has a cellular
structure for each manifold and one uses the map to $X$ to label
cells. State sums are then defined (over appropriate labellings) for
closed manifolds and then extended to homotopy quantum field
theories. Rodrigues views this as a way of incorporating matter into
lattice topological quantum field theory.

\subsection*{Examples from variants of modular categories} Turaev also
constructs 2+1-dimensional theories from algebraic data \cite{tu:hft}
by defining the notion of a modular $\pi$-category which is a
generalisation of a modular category and when $X$ is contractible and
the construction reduces to his construction of a 2+1-dimensional TQFT
from a modular category.

\section*{Thin invariant field theories}
There is already available a notion of ``higher'' line bundle with
connection in the form of a {\em gerbe} and in this section we explain
the work of U. Bunke and the second two authors \cite{butuwi:ghqft}
relating these to homotopy quantum field theories.

Following Barrett's work on classifying vector bundles by holonomy
\cite{Barrett:HolonomyAndPathStructures} one can modify the definition
of the homotopy surface category of a smooth finite dimensional
manifold by considering collared smooth cobordisms and demanding
invariance under smooth homotopies with rank $\leq 2$. Such homotopies
are known as {\em thin homotopies} and may be informally characterised
as those sweeping out no volume. Call this category the {\em
  thin-homotopy surface category} $\cT_X$ and define a {\em rank one
  thin invariant field theory with target a smooth manifold $X$} to be
a symmetric multiplicative functor $E$ from $\cT_X$ to the category of
one-dimensional vector spaces such that there exists a closed three
from $c$ on $X$ satisfying $E (\partial v) = \exp(i{\textstyle\int_V
  v^*c})$ whenever $v\colon V\ra X$ is an $X$-three-manifold. This
last condition guarantees that the functor $E$ is suitably smooth.  If
the three form $c$ is zero we say the thin invariant field theory is
{\em flat}. Though it is not obvious, a flat thin invariant field
theory is the same thing as a rank one, normalised homotopy quantum
field theory. Here {\em normalised} means that for each point, the
induced TQFT is trivial.

Examples   of thin invariant field  theories   can be constructed from
gerbes with connection. A  gerbe with connection  may be  described by
its         holonomy                 (see           for        example
\cite{MackaayPicken:HolonomyAndParallelTransport}) , which in turn can
be viewed as a  $\bC^\times$-valued function $S$ on  the space of maps
of closed surfaces into $X$. Given such,  and supposing that $H_1X$ is
trivial, it is possible to construct a \tift over $X$ as follows.
For a loop $\gamma$ set
\[
E(\gamma) = \mbox{Span}_\bC \{\sigma\colon W \ra X \; | \;
\partial\sigma = \gamma  \} / \sim
\]
where $\sigma \sim \sigma^\prime$ if $\sigma = S(\langle \sigma \cup
{\sigma^\prime}^{op}\rangle )\sigma^\prime$, where the notation
$\langle -\rangle$ indicates that we are viewing 
$\sigma \cup {\sigma^\prime}^{op} $ as a closed surface. For a
cobordism $g\colon \Sigma\ra X$ from $\gamma$ to $\gamma^\prime$ the
map $E(\gamma) \ra E(\gamma^\prime)$ is given by composing surfaces as
indicated below, where $\sigma\colon W \ra X$ is a generator of
$E(\gamma)$.

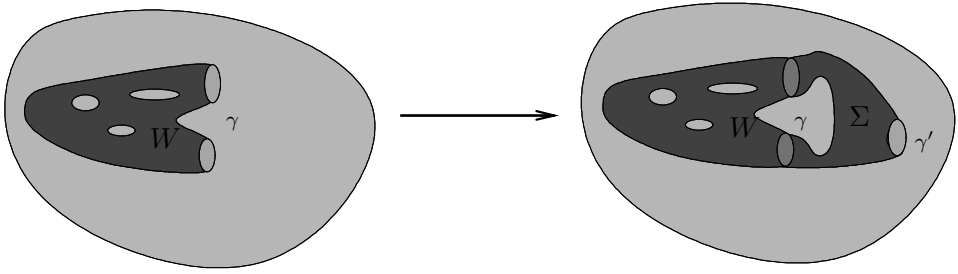
\vspace{0.3cm}
\begin{center}
\input{comp2.eepic}
\end{center}
\vspace{0.3cm}

In fact , for a background space satisfying $H_1X=0$, this
construction sets up a one-to-one correspondence between gerbes with
connection on $X$ and \tifts on $X$.

By virtue of the remarks above about flat \tifts and rank one
normalised homotopy quantum field theories, another source of (flat)
\tifts is Turaev's construction using cocycles.

With some more work it can be shown (see \cite{butuwi:ghqft}) that
for a general finite dimensional smooth manifold $X$ there is a
correspondence 

\[
\{\mbox{Gerbes with connection on $X$}\} \leftrightarrow 
\{\mbox{ \tifts on $X$}\} 
\]

In other words a \tift can be viewed as an alternative
characterisation of a gerbe. Flat gerbes, it turns out, correspond to
flat \tifts and hence to rank one, normalised HQFTs. By identifying
flat gerbes with $H^2(X;\bC^\times)$, this identifies rank one,
normalised HQFTs on $X$ with $H^2(X;\bC^\times)$. This extends a 
result of Turaev who showed this for $X$ and Eilenberg-Maclane space
using the algebraic classification in terms of crossed algebras.

\section*{Extended \hqfts}
In this section we will discuss another variant of homotopy quantum
field theory. Elsewhere in this volume Stephen Sawin discusses
2-TQFTs, which extend TQFTs to a three tier structure rather than the
usual two tier structure. A version of \hqft along similar lines has
been defined by the first two authors in
\cite{BrightwellTurner:enriched}.  The approach is to begin by
defining the 2-category version of the homotopy surface category of a
space and then considering some kind of representations of this. This
approach has been pioneered in the background-free case by Tillmann in
\cite{Tillmann:Discrete, ti:sxs}.  She considers a 2-category whose
objects are one-manifolds, whose morphisms are surfaces and whose
2-morphisms are diffeomorphisms. She explains how this may be
thought of as a discrete approximation to Segal's topological category
occurring in the definition of CFT, which has Riemann
surfaces as morphisms. Indeed, the classifying space of a
morphism category in Tillmann's 2-category has the same rational
homotopy type as the morphism space in Segal's category, thus
making the link between TQFTs extended to the 2-category setting as
above and CFT. Generalising this to \hqfts is then some way toward
approximating the elliptic objects discussed in the first section.
The {\em homotopy surface 2-category of $X$} is the 2-category
$\cS_X^{(2)}$ defined
(roughly) as follows. 
\begin{itemize}
\item Objects are based continuous functions $s\colon S_m \rightarrow X$,
  for $m\in \bN$. 
\item 1-Morphisms are  continuous functions $g\colon \Sigma \rightarrow
  X$ where $\Sigma$ is a surface and $g$ agrees with source and target
  boundaries.
\item 2-Morphisms are orientation preserving diffeomorphisms $T\colon
  \Sigma_1 \rightarrow \Sigma_2$ that fix boundary collars pointwise
  and such that the following diagram commutes up to
  basepoint-preserving homotopy relative to the boundary.
\[
\xymatrix{\Sigma_1 \ar[rr]^{T} \ar[dr]_{g_1} & & \Sigma_2
  \ar[dl]^{g_2}\\ & X}
\]
\end{itemize}

2-morphisms are identified if they are in the same component of the
mapping class group and 1-morphisms are identified up to limited
isotopy. It would be inappropriate in this short review to present
more details (the 2-category-minded reader will be aware there must be
many) and the reader is referred to \cite{BrightwellTurner:enriched}.  
An {\em extended \hqft}
is then defined to be a multiplicative strict 2-functor from $\cS_X^{(2)}$
to the 2-category of additive, idempotent complete $\bC$-linear
categories. There are other possibilites for the target two category such
as 2-vector spaces or those considered by Sawin.

Suppose that $X$ is simply connected and let $\cA$ be the category
associated to the circle equipped with the constant map to some
basepoint in $X$. This category inherits a rich algebraic structure,
in particular it is a semi-simple balanced category (monoidal with
braiding and twist) with an action of $\pi_2(X)$ satisfying a number
of further conditions describing the interaction of the action with
the balanced structure. The category $\cA$ also has defined on it an
involution $(-)^*\colon \cA \ra \cA$ (defined by Tillmann) and it is
shown in \cite{BrightwellTurner:enriched} that for certain self-dual
extended \hqfts this provides $\cA$ with a right duality, thus turning
$\cA$ into a lax semi-simple tortile category with $\pi_2(X)$-action.

There is a similar situation for the case when $X$ is an
Eilenberg-Maclane space of type $K(\pi,1)$, where a self-dual extended
\hqft gives rise to a lax version of a tortile $\pi$-category. A
strict tortile $\pi$-category is the same thing as a ribbon
$\pi$-category defined by Turaev, showing that  self dual
1+1-dimensional extended \hqfts have the expected close connection
with and 2+1-dimensional homotopy quantum field theories.

\end{document}